\documentclass[12pt]{article}
\usepackage{latexsym,amsfonts,amsmath,amssymb}
\newtheorem{theorem}{Theorem}
\newtheorem{corollary}[theorem]{Corollary}
\newtheorem{sublemma}{Lemma}[theorem]
\newtheorem{lemma}[theorem]{Lemma}
\newtheorem{question}[theorem]{Question}
\newtheorem{observation}[theorem]{Observation}
\newtheorem{claim}[theorem]{Claim}
\newtheorem{subclaim}{Claim}[sublemma]
\newtheorem{conjecture}[theorem]{Conjecture}
\newtheorem{fact}[theorem]{Fact}
\newtheorem{definition}[theorem]{Definition}
\newtheorem{remark}[theorem]{Remark}
\newtheorem{example}[theorem]{Example}
\newtheorem{exercise}{Exercise}[section]
\def\Theorem #1.#2 #3\par{\setbox1=\hbox{#1}\ifdim\wd1=0pt
   \begin{theorem}{\rm #2} #3\end{theorem}\else
   \newtheorem{#1}[theorem]{#1}\begin{#1}\label{#1}{\rm #2} #3\end{#1}\fi}
\def\Corollary #1.#2 #3\par{\setbox1=\hbox{#1}\ifdim\wd1=0pt
   \begin{corollary}{\rm #2} #3\end{corollary}\else
   \newtheorem{#1}[theorem]{#1}\begin{#1}\label{#1}{\rm #2} #3\end{#1}\fi}
\def\Lemma #1.#2 #3\par{\setbox1=\hbox{#1}\ifdim\wd1=0pt
   \begin{lemma}{\rm #2} #3\end{lemma}\else
   \newtheorem{#1}[theorem]{#1}\begin{#1}\label{#1}{\rm #2} #3\end{#1}\fi}
\def\SubLemma #1.#2 #3\par{\setbox1=\hbox{#1}\ifdim\wd1=0pt
   \begin{sublemma}{\rm #2} #3\end{sublemma}\else
   \newtheorem{#1}{#1}[theorem]\begin{#1}\label{#1}{\rm #2} #3\end{#1}\fi}
\def\Question #1.#2 #3\par{\setbox1=\hbox{#1}\ifdim\wd1=0pt
   \begin{question}{\rm #2} #3\end{question}\else
   \newtheorem{#1}[theorem]{#1}\begin{#1}\label{#1}{\rm #2} #3\end{#1}\fi}
\def\Observation #1.#2 #3\par{\setbox1=\hbox{#1}\ifdim\wd1=0pt
   \begin{observation}{\rm #2} #3\end{observation}\else
   \newtheorem{#1}[theorem]{#1}\begin{#1}\label{#1}{\rm #2} #3\end{#1}\fi}
\def\Claim #1.#2 #3\par{\setbox1=\hbox{#1}\ifdim\wd1=0pt
   \begin{claim}{\rm #2} #3\end{claim}\else
   \newtheorem{#1}[theorem]{#1}\begin{#1}\label{#1}{\rm #2} #3\end{#1}\fi}
\def\SubClaim #1.#2 #3\par{\setbox1=\hbox{#1}\ifdim\wd1=0pt
   \begin{subclaim}{\rm #2} #3\end{subclaim}\else
   \newtheorem{#1}{#1}[sublemma]\begin{#1}\label{#1}{\rm #2} #3\end{#1}\fi}
\def\Conjecture #1.#2 #3\par{\setbox1=\hbox{#1}\ifdim\wd1=0pt
   \begin{conjecture}{\rm #2} #3\end{conjecture}\else
   \newtheorem{#1}[theorem]{#1}\begin{#1}\label{#1}{\rm #2} #3\end{#1}\fi}
\def\Fact #1.#2 #3\par{\setbox1=\hbox{#1}\ifdim\wd1=0pt
   \begin{fact}{\rm #2} #3\end{fact}\else
   \newtheorem{#1}[theorem]{#1}\begin{#1}\label{#1}{\rm #2} #3\end{#1}\fi}
\def\Definition #1.#2 #3\par{\setbox1=\hbox{#1}\ifdim\wd1=0pt
   \begin{definition}{\rm #2} {\rm #3}\end{definition}\else
   \newtheorem{#1}[theorem]{#1}\begin{#1}\label{#1}{\rm #2} {\rm #3}\end{#1}\fi}
\def\Remark #1.#2 #3\par{\setbox1=\hbox{#1}\ifdim\wd1=0pt
   \begin{remark}{\rm #2} {\rm #3}\end{remark}\else
   \newtheorem{#1}[theorem]{#1}\begin{#1}\label{#1}{\rm #2} {\rm #3}\end{#1}\fi}
\def\Example #1.#2 #3\par{\setbox1=\hbox{#1}\ifdim\wd1=0pt
   \begin{example}{\rm #2} #3\end{example}\else
   \newtheorem{#1}[theorem]{#1}\begin{#1}\label{#1}{\rm #2} #3\end{#1}\fi}
\def\Exercise #1.#2 #3\par{\setbox1=\hbox{#1}\ifdim\wd1=0pt
   {\footnotesize\begin{exercise}{\rm #2} {\rm #3}\end{exercise}}\else
   \newtheorem{#1}[section]{#1}{\footnotesize\begin{#1}\label{#1}{\rm #2} {\rm #3}\end{#1}}\fi}
\def\QuietTheorem #1.#2 #3\par{\setbox1=\hbox{#1}\ifdim\wd1=0pt\proclaim{Theorem {\rm #2}}{#3}\else\proclaim{#1 {\rm #2}}{#3}\fi}
\newcommand{\proclaim}[2]{\smallskip\noindent{\bf #1} {\sl#2}\par\smallskip}
\def\Proclaim #1.#2 #3\par{\proclaim{#1 {\rm #2}}{#3}}
\newenvironment{proof}{\noindent}{\kern2pt\QEDbox\par\bigskip}
\def\Proof#1: {\setbox1=\hbox{#1}\ifdim\wd1=0pt\begin{proof}{\bf Proof: }\else\medskip\begin{proof}{\bf #1: }\fi}
\newcommand{\QED}{\end{proof}}
\def\BF#1.{{\bf #1.}}
\def\Abstract #1\par{\begin{quotation}{\singlespaced\footnotesize{\noindent{\bf Abstract.~}#1}}\end{quotation}}
\def\Title #1\par{\title{#1}\maketitle}
\def\Author #1\par{\author{#1}}
\def\Acknowledgement#1\par{\thanks{#1}}
\def\Chapter #1\par{\chapter{#1}}
\def\Section #1\par{\section{#1}}
\def\QuietSection #1\par{\section*{#1}}
\def\SubSection #1\par{\subsection{#1}}
\def\SubSubSection #1\par{\subsubsection{#1}}
\def\MidTitle #1\par{\bigskip\goodbreak\centerline{\small\bf #1}\bigskip\noindent}
\def\Margin #1\par{\marginpar{\tiny #1}}

\newcommand{\singlespaced}{\baselineskip=15pt}
\def\bottomnote #1\par{{\renewcommand{\thefootnote}{}\footnotetext{#1}}}
\newcommand{\A}{{\mathbb A}}
\newcommand{\B}{{\mathbb B}}
\newcommand{\C}{{\mathbb C}}

\renewcommand{\P}{{\mathbb P}}
\newcommand{\Q}{{\mathbb Q}}
\def\R{{\mathbb R}}

\newcommand{\term}{{\!\scriptscriptstyle\rm term}}
\newcommand{\Dterm}{{D_{\!\scriptscriptstyle\rm term}}}

\newcommand{\Hterm}{{H_{\!\scriptscriptstyle\rm term}}}

\newcommand{\Qterm}{{\Q_{\!\scriptscriptstyle\rm term}}}
\newcommand{\Gterm}{{G_{\!\scriptscriptstyle\rm term}}}

\newcommand{\Vbar}{{\overline{V}}}

\renewcommand{\Ddot}{{\dot D}}
\newcommand{\Pdot}{{\dot\P}}
\newcommand{\Qdot}{{\dot\Q}}
\newcommand{\qdot}{{\dot q}}

\newcommand{\Rdot}{{\dot\R}}

\newcommand{\xdot}{{\dot x}}

\newcommand{\one}{\mathop{1\hskip-3pt {\rm l}}}
\newfont{\msam}{msam10 at 12pt}

\newcommand{\of}{\subseteq}

\newcommand{\set}[1]{\{\,{#1}\,\}}

\newcommand{\compose}{\circ}

\newcommand{\elesub}{\prec}

\newcommand{\dom}{\mathop{\rm dom}}

\newcommand{\add}{\mathop{\rm add}}
\newcommand{\coll}{\mathop{\rm coll}}
\newcommand{\cof}{\mathop{\rm cof}}

\newcommand{\con}{\mathop{\hbox{\sc con}}}

\newcommand{\plus}{{+}}

\newcommand{\restrict}{\upharpoonright}
\newcommand{\satisfies}{\models}
\newcommand{\forces}{\Vdash}

\newcommand{\necessary}{\mathop{\raisebox{-1pt}{$\Box$\!}}}
\newcommand{\cross}{\times}

\newcommand{\union}{\cup}
\newcommand{\Union}{\bigcup}
\newcommand{\intersect}{\cap}

\newcommand{\smalllt}{\mathrel{\mathchoice{\raise2pt\hbox{$\scriptstyle<$}}{\raise1pt\hbox{$\scriptstyle<$}}{\scriptscriptstyle<}{\scriptscriptstyle<}}}
\newcommand{\smallleq}{\mathrel{\mathchoice{\raise2pt\hbox{$\scriptstyle\leq$}}{\raise1pt\hbox{$\scriptstyle\leq$}}{\scriptscriptstyle\leq}{\scriptscriptstyle\leq}}}
\newcommand{\ltomega}{{{\smalllt}\omega}}

\newcommand{\ltkappa}{{{\smalllt}\kappa}}

\newcommand{\boolval}[1]{\mathopen{\lbrack\!\lbrack}\,#1\,\mathclose{\rbrack\!\rbrack}}
\def\[#1]{\boolval{#1}}

\newcommand{\UnderTilde}[1]{{\setbox1=\hbox{$#1$}\baselineskip=0pt\vtop{\hbox{$#1$}\hbox to\wd1{\hfil$\sim$\hfil}}}{}}
\newcommand{\Undertilde}[1]{{\setbox1=\hbox{$#1$}\baselineskip=0pt\vtop{\hbox{$#1$}\hbox to\wd1{\hfil$\scriptstyle\sim$\hfil}}}{}}
\newcommand{\undertilde}[1]{{\setbox1=\hbox{$#1$}\baselineskip=0pt\vtop{\hbox{$#1$}\hbox to\wd1{\hfil$\scriptscriptstyle\sim$\hfil}}}{}}
\newcommand{\UnderdTilde}[1]{{\setbox1=\hbox{$#1$}\baselineskip=0pt\vtop{\hbox{$#1$}\hbox to\wd1{\hfil$\approx$\hfil}}}{}}
\newcommand{\Underdtilde}[1]{{\setbox1=\hbox{$#1$}\baselineskip=0pt\vtop{\hbox{$#1$}\hbox to\wd1{\hfil\scriptsize$\approx$\hfil}}}{}}

\newcommand{\st}{\mid}
\renewcommand{\th}{{\hbox{\scriptsize th}}}
\newcommand{\Iff}{\mathrel{\leftrightarrow}}

\newcommand{\iso}{\cong}
\def\<#1>{\langle#1\rangle}

\newcommand{\QEDbox}{\fbox{}}

\newcommand{\ORD}{\mathop{\hbox{\sc ord}}}

\newcommand{\ZFC}{\hbox{\sc zfc}}

\newcommand{\AD}{\hbox{\sc ad}}

\newcommand{\MA}{\hbox{\sc ma}}

\newcommand{\MP}{\hbox{\sc mp}}

\newcommand{\MPtilde}{\UnderTilde{\MP}}
\newcommand{\MPccc}{{\MP_{\scriptsize\!\hbox{\sc ccc}}}}

\newcommand{\cell}[1]{\boxit{\hbox to 17pt{\strut\hfil$#1$\hfil}}}
\newcommand{\head}[2]{\lower2pt\vbox{\hbox{\strut\footnotesize\it\hskip3pt#2}\boxit{\cell#1}}}
\newcommand{\boxit}[1]{\setbox4=\hbox{\kern2pt#1\kern2pt}\hbox{\vrule\vbox{\hrule\kern2pt\box4\kern2pt\hrule}\vrule}}
\newcommand{\Col}[3]{\hbox{\vbox{\baselineskip=0pt\parskip=0pt\cell#1\cell#2\cell#3}}}
\newcommand{\tapenames}{\raise 5pt\vbox to .7in{\hbox to .8in{\it\hfill input: \strut}\vfill\hbox to
.8in{\it\hfill scratch: \strut}\vfill\hbox to .8in{\it\hfill output: \strut}}}
\newcommand{\Head}[4]{\lower2pt\vbox{\hbox to25pt{\strut\footnotesize\it\hfill#4\hfill}\boxit{\Col#1#2#3}}}
\newcommand{\Dots}{\raise 5pt\vbox to .7in{\hbox{\ $\cdots$\strut}\vfill\hbox{\ $\cdots$\strut}\vfill\hbox{\
$\cdots$\strut}}}
\renewcommand{\dots}{\raise5pt\hbox{\ $\cdots$}}
\newcommand{\factordiagramup}[6]{$$\begin{array}{ccc}
#1&\raise3pt\vbox{\hbox to60pt{\hfill$\scriptstyle
#2$\hfill}\vskip-6pt\hbox{$\vector(4,0){60}$}}&#3\\ \vbox
to30pt{}&\raise22pt\vtop{\hbox{$\vector(4,-3){60}$}\vskip-22pt\hbox
to60pt{\hfill$\scriptstyle #4\qquad$\hfill}}
     &\ \ \lower22pt\hbox{$\vector(0,3){45}$}\ {\scriptstyle #5}\\
\vbox to15pt{}&&#6\\
\end{array}$$}
\newcommand{\factordiagram}[6]{$$\begin{array}{ccc}
#1&&\\ \ \ \raise22pt\hbox{$\vector(0,-3){45}$}\ {\scriptstyle #2}
&\raise22pt\hbox{$\vector(2,-1){90}$}\raise5pt\llap{$\scriptstyle#3$\qquad\quad}&\vbox
to25pt{}\\ #4&\raise3pt\vbox{\hbox to90pt{\hfill$\scriptstyle
#5$\hfill}\vskip-6pt\hbox{$\vector(4,0){90}$}}&#6\\
\end{array}$$}
\newcommand{\df}{\it} 
\begin{document}
\author{Joel David Hamkins\\
\normalsize\sc The City University of New York\\
{\footnotesize http://jdh.hamkins.org}\\
\\
W. Hugh Woodin\\
\normalsize\sc University of California at Berkeley\\
{\footnotesize http://math.berkeley.edu/$\sim$woodin}} \bottomnote MSC: 03E55, 03E40. Keywords: forcing axiom, ccc forcing, weakly compact cardinal.
The first author is affiliated with the College of Staten Island of CUNY and The CUNY Graduate Center, and his research has been supported by grants
from the Research Foundation of CUNY and the National Science Foundation DMS-9970993. The research of the second author has been partially supported
by National Science Foundation grant DMS-9970255.

\Title The Necessary Maximality Principle for c.c.c.~forcing is equiconsistent with a weakly compact cardinal

\Abstract The Necessary Maximality Principle for c.c.c.~forcing with real parameters is equiconsistent with the existence of a weakly compact
cardinal.

The Necessary Maximality Principle for c.c.c.~forcing, denoted $\necessary\MPccc(\R)$, asserts that any statement about a real in a c.c.c.~extension
that could become true in a further c.c.c.~extension and remain true in all subsequent c.c.c.~extensions, is already true in the minimal extension
containing the real. We show that this principle is equiconsistent with the existence of a weakly compact cardinal.

The principle is one of a family of principles considered in \cite{Hamkins2003:MaximalityPrinciple} (building on ideas of
\cite{Chalons:FullSetTheory} and overlapping with independent work in \cite{StaviVaananen:ReflectionPrinciples}). The family begins with the
Maximality Principle $\MP$, the scheme asserting the truth of any statement that holds in some forcing extension $V^\P$ and all subsequent extensions
$V^{\P*\Qdot}$ (these are the {\df forceably necessary} statements).  The boldface form $\MPtilde$ allows real parameters in the scheme, and the
Necessary Maximality Principle $\necessary\MPtilde$ asserts $\MPtilde$ in all forcing extensions, using the parameters available in those extensions.
The main results of \cite{Hamkins2003:MaximalityPrinciple} show that $\MP$ is equiconsistent with \ZFC, while $\MPtilde$ is equiconsistent with the
L\'evy scheme ``$\ORD$ is Mahlo'' and $\necessary\MPtilde$ is far stronger. Philip Welch proved that $\necessary\MPtilde$ implies Projective
Determinacy, and the second author of this paper improved the conclusion to $\AD^{L(\R)}$. He also provided an upper bound by proving the consistency
of $\necessary\MPtilde$ from the theory ``$\AD_\R+\Theta$ is regular''.

In this article, we focus on the principles obtained by restricting attention to the class of c.c.c.~forcing notions. The parameter-free version
$\MPccc$ asserts the truth of any statement holding in some c.c.c.~extension $V^\P$ and all subsequent c.c.c.~extensions $V^{\P*\Qdot}$. This is
equiconsistent with \ZFC\ by \cite[Corollary 32]{Hamkins2003:MaximalityPrinciple}. An almost identical principle, where one requires the upward
absoluteness of the statement from $V^\P$ to any $V^{\P*\Qdot}$ to be \ZFC-provable (rather than merely true), was considered independently in
\cite{StaviVaananen:ReflectionPrinciples}.

The principle $\MPccc$ implies a spectacular failure of the Continuum Hypothesis. The reason is that with c.c.c.~forcing one may add as many Cohen
reals as desired, and once they are added, of course, the value of the continuum $2^\omega$ remains inflated in all subsequent c.c.c.~extensions.
Thus, the assertion that $2^\omega$ is larger than $\aleph_{\omega^{17}}$, say, or any cardinal whose definition is c.c.c.~absolute, is
c.c.c.~forceably necessary, and hence true under $\MPccc$.

For the boldface version of $\MPccc$, there is initially little reason to restrict as in $\MPtilde$ to real parameters, and so we denote by
$\MPccc(X)$ the scheme in which arbitrary parameters in $X$ are allowed. Because for any parameter $z$ the assertion $|z|<2^\omega$ is
c.c.c.~forceably necessary, we can't allow parameters outside $H(2^\omega)$. Parameters inside $H(2^\omega)$, however, are fine, and
$\MPccc(H(2^\omega))$ is equiconsistent with $\MPtilde$, which as we have mentioned is equiconsistent with the L\'evy scheme (see
\cite{Hamkins2003:MaximalityPrinciple}). The weaker principle $\MPccc(\R)$ has recently been proved by Leibman \cite{Leibman2004:Dissertation} to be
equiconsistent with \ZFC, and one may freely add a large initial segment of the ordinals as parameters.

The strongest form of the principle is $\necessary\MPccc(X)$, which asserts that $\MPccc(X)$ holds in all c.c.c.~extensions, reinterpreting $X$ {\it
de dicto} in these extensions. Thus, the principle asserts that if $x$ is in $X$ in some c.c.c.~extension $V^{\P_0}$ and $\varphi(x)$ holds in a
further c.c.c.~extension $V^{\P_0*\P}$ and all subsequent c.c.c.~extensions $V^{\P_0*\P*\Qdot}$, then $\varphi(x)$ holds already in $V^{\P_0}$.
Because $\MPccc(H(2^\omega))$ is equiconsistent with $\MPtilde$, one might have expected the same for $\necessary\MPccc(H(2^\omega))$ and
$\necessary\MPtilde$. But the former principle is simply false.

\Observation.({Leibman \cite{Leibman2004:Dissertation}}) $\necessary\MPccc(H(2^\omega))$  is false.\label{GeorgeObservation}

Leibman merely observed that $\MPccc(H(2^\omega))$ implies Martin's Axiom \MA, because the assertion that there is a filter for a given
c.c.c.~partial order meeting a certain family of dense sets is c.c.c.~forceably necessary. Thus, if $\necessary\MPccc(H(2^\omega))$ held, then \MA\
would hold in all c.c.c.~extensions. But \MA\ does not hold in all c.c.c.~extensions, because even the forcing to add a single Cohen real creates
Souslin trees. This argument makes an essential use of the uncountable parameters available in $H(2^\omega)$, such as the Souslin trees in the Cohen
extension, and there appears to be no general way to get by with just real parameters (although doing so in the special case when $\omega_1$ is
accessible to reals is the key to Theorem \ref{NMPcccOmega1Inacc} below).

So when it comes to the necessary form of the principle, the natural collection of parameters is $\R$ after all, and we focus our attention on the
principle $\necessary\MPccc(\R)$.

\Question Main Question. Is $\necessary\MPccc(\R)$ consistent?

This question is answered by our main theorem.


\Theorem Main Theorem. The principle $\necessary\MPccc(\R)$ is equiconsistent over \ZFC\ with the existence of a weakly compact cardinal.

The rest of this article consists of our proof of this theorem, followed by a short application of the proof to $\MP(\R)$. We concentrate first on
the converse direction of the Main Theorem. Let $V_\delta\elesub V$ denote the scheme, in the language with a constant symbol for $\delta$, asserting
for every formula $\varphi$ in the language of set theory that $\forall x\in V_\delta\,[\varphi(x)\Iff\varphi(x)^{V_\delta}]$. The point of this is
that in the construction of Theorem \ref{ForcingNMPccc} we would like at heart to have a truth predicate for $V$, which is of course lacking by
Tarski's theorem, but we can get by merely with a truth predicate for $V_\delta$ and the scheme $V_\delta\elesub V$. Note that if $V_\delta\elesub V$
and $G\of\P$ is $V$-generic for forcing $\P\in V_\delta$, then $V_\delta[G]\elesub V[G]$, because $V_\delta$ and $V$ agree on whether a given
statement is forced.

\Lemma. If there is a model of $\ZFC+{}$there is a weakly compact cardinal, then there is a model of $\ZFC+{}$there is a weakly compact
cardinal${}+V_\delta\elesub V$.\label{Reflection}

\Proof: Let $T$ be the latter theory, and suppose that $M$ is a model of \ZFC\ with a weakly compact cardinal. Since $M$ satisfies every instance of
the L\'evy Reflection Theorem, it follows that every finite subset of $T$ is consistent, by interpreting $\delta$ to be a sufficiently reflective
ordinal of $M$. And so $T$ as a whole is consistent.\QED

The converse implication of the Main Theorem now follows from:

\Theorem. Assume $\kappa$ is weakly compact, $\kappa<\delta$ and $V_\delta\elesub V$. Then there is a forcing extension satisfying
$\necessary\MPccc(\R)+\kappa=\omega_1$.\label{ForcingNMPccc}

The proof of this theorem relies in part on some general facts due to Kunen and Harrington-Shelah
\cite{HarringtonShelah1985:SomeExactEquiconsistencyResults} concerning forcing and weakly compact cardinals. For completeness, we include proofs
here.

\SubLemma. If $\kappa$ is weakly compact, then any finite support product of $\kappa$-c.c.~forcing is $\kappa$-c.c.\label{ProductIsKappaCC}

\Proof: Suppose first that $\P$ and $\Q$ are $\kappa$-c.c.~and that $A\of\P\cross\Q$ is an antichain of size $\kappa$ in the product. Enumerate
$A=\set{(p_\alpha,q_\alpha)\st\alpha<\kappa}$ and define the coloring $f:[\kappa]^2\to 2$ by $f(\alpha,\beta)=0$ if $p_\alpha\perp p_\beta$,
otherwise $1$. Since $\kappa$ is weakly compact, there is a homogeneous set $H\of\kappa$ of size $\kappa$, meaning that $f$ is constant on $[H]^2$.
If the constant value is $0$, then $p_\alpha\perp p_\beta$ for all $\alpha,\beta\in H$, contradicting that $\P$ is $\kappa$-c.c. Otherwise the
constant value is $1$, in which case $q_\alpha\perp q_\beta$ for all such $\alpha$ and $\beta$, contradicting that $\Q$ is $\kappa$-c.c. By
induction, it follows that any finite product of $\kappa$-c.c.~forcing is $\kappa$-c.c. Consider now an antichain $A$ of size $\kappa$ in an
arbitrary finite-support product $\Pi_{\alpha\in I}\P_\alpha$, where each $\P_\alpha$ is $\kappa$-c.c. By a delta system argument, we may assume that
supports of the conditions in $A$ form a delta system. Any two conditions in $A$ must be incompatible on the root of this system, contradicting the
fact that any finite product of $\kappa$-c.c.~partial orders is $\kappa$-c.c.\QED

\SubLemma. If $\kappa$ is weakly compact and $\B$ is a $\kappa$-c.c.~complete Boolean algebra, then every subset $A\of\B$ of size less than $\kappa$
generates a complete subalgebra that is also of size less than $\kappa$.

\Proof: Construct the increasing continuous sequence of subalgebras $A=A_0\of A_1\of\cdots\of A_\alpha\of\cdots$, for $\alpha<\kappa$, where
$A_{\alpha+1}$ is obtained by adding to $A_\alpha$ the infima (computed in $\B$) of all subsets of $A_\alpha$ and closing under the Boolean
operations, taking unions at limits. Let $\A=\union_\alpha A_\alpha$. Since $\B$ is $\kappa$-c.c., all the antichains of $\A$ live on some
$A_\alpha$, and so $\A$ is the complete subalgebra generated by $A$ in $\B$. Since $|A_\alpha|<\kappa$ for all $\alpha<\kappa$, it follows that
$|\A|\leq\kappa$. By moving to an isomorphic copy of $\B$, assume $\A\of\kappa$ and place $\A$ into a transitive model of set theory $M$ of size
$\kappa$. Fix a weakly compact embedding $j:M\to N$ with critical point $\kappa$, and observe in $N$ that $\A$ is a complete subalgebra of $j(\A)$
containing the generators $A$. Thus, $\A=j(\A)$ and so $|\A|<\kappa$, as desired.\QED

\SubLemma. Consequently, if $\kappa$ is weakly compact and $G\of\P$ is $V$-generic for $\kappa$-c.c.~forcing $\P$, then every $x\in H(\kappa)^{V[G]}$
is $V$-generic for $\kappa$-c.c.~forcing of size less than $\kappa$.\label{SmallSuborder}

\Proof: Let $\B$ be the regular open algebra of $\P$. By coding, we may assume that $x\of\beta$ for some $\beta<\kappa$. Let $\xdot$ be a $\P$-name
for $x$ such that $\forces_\P\xdot\of\check\beta$. For each $\xi<\beta$, let $b_\xi=\boolval{\check\xi\in\xdot}$. By the previous lemma, the complete
subalgebra $\A$ generated by $A=\set{b_\xi\st\xi<\beta}$ has size less than $\kappa$. And clearly $x$ is constructible from $G\intersect A$.\QED

Our proof also relies on the term forcing construction, which we now explain. Suppose that $\P$ is any partial order and $\Qdot$ is the $\P$-name of
a partial order. The term forcing $\Qterm$ for $\Qdot$ over $\P$ consists of conditions $q$ such that $\one\forces q\in\Qdot$, with the order
$p\leq_\term q$ if and only if $\one\forces p\leq_\Qdot q$. One can restrict the size of $\Qterm$ by using only the names $p$ in a full set $B$ of
names, meaning that for any $\P$-name $q$ with $\one\forces q\in\Qdot$ there is $p\in B$ with $\one\forces p=q$. Any such full set of names forms a
dense subset of $\Qterm$, which is therefore equivalent as a forcing notion. The fundamental property of term forcing is the following:

\SubLemma. Suppose that $\Hterm\of\Qterm$ is $V$-generic for the term forcing of $\Qdot$ over $\P$ and $\Vbar$ is any model of set theory with
$V[\Hterm]\of\Vbar$. If there is a $V$-generic filter $G\of\P$ in $\Vbar$, then there is a $V[G]$-generic filter $H\of\Q=\Qdot_G$ in
$\Vbar$.\label{TermForcing}

\Proof: The filter $\Hterm$ consists of $\P$-names for conditions in $\Qdot$, so it makes sense to let $H=\set{q_G\st q\in\Hterm}$ in $\Vbar$. To see
that this is $V[G]$-generic for $\Q=\Qdot_G$, suppose that $D\of\Q$ is a dense subset of $\Q$ in $V[G]$. Let $\Ddot$ be a $\P$-name for $D$, forced
by $\one$ to be dense, and let $\Dterm$ be the set of conditions $\qdot$ that are forced by $\one$ to be in $\Ddot$. It is easy to see that $\Dterm$
is a dense subset of $\Qterm$, and so there is a condition $q\in\Hterm\intersect\Dterm$. It follows that $q_G\in H\intersect D$, and so $H$ meets
$D$, as desired.\QED

In particular, if $\Vbar=V[\Hterm][B]$ is a forcing extension of $V[\Hterm]$ containing a $V$-generic $G\of\P$, then we may rearrange the forcing as
$\Vbar=V[\Hterm][B]=V[G][H][(\Hterm*B)/(G*H)]$, where $(\Hterm*B)/(G*H)$ is the quotient forcing adding the rest of $\Hterm*B$ over $V[G][H]$. There
is no need in Lemma \ref{TermForcing} for $G$ and $\Hterm$ to be mutually $V$-generic.

\SubLemma. In the context of the previous Lemma, if $\kappa$ is weakly compact in $V$, $|\P|<\kappa$ and $\one\forces_\P\Qdot$ is
$\check\kappa$-c.c., then $\Qterm$ is $\kappa$-c.c.~in $V$ (and hence also in $V[G]$). If $\Vbar=V[\Hterm][B]$ is a $\kappa$-c.c.~forcing extension
of $V[\Hterm]$, then it is a $\kappa$-c.c.~extension of the resulting $V[G][H]$.\label{TermForcingKappaCC}

\Proof: Suppose that $A\of\Qterm$ is an antichain in $V$ of size $\kappa$. Any two elements $q_0,q_1\in A$ are incompatible in $\Qterm$, meaning that
there is no condition $q\in\Qterm$ such that $\one_\P$ forces both $q\leq_\Qdot q_0$ and $q\leq_\Qdot q_1$. It follows that there is some condition
$p\in\P$ such that $p\forces_\P q_0\perp_\Qdot q_1$. Enumerate $A=\set{q_\beta\st\beta<\kappa}$ and define $f:[\kappa]^2\to\P$ by $f(\alpha,\beta)=p$
for some $p$ forcing $q_\alpha\perp_\Qdot q_\beta$. Since $\kappa$ is weakly compact, there is a homogeneous set $H\of\kappa$ of size $\kappa$ on
which $f$ has some constant value $p$. Thus, $p\forces q_\alpha\perp_\Qdot q_\beta$ for all $\alpha<\beta$ from $H$, contradicting that $\one$ forces
$\Qdot$ is $\check\kappa$-c.c. So $\Qterm$ is $\kappa$-c.c.~in $V$.

If $\Vbar=V[\Hterm][B]$ is a forcing extension of $V[\Hterm]$, then we have already observed that $\Vbar=V[G][H][(\Hterm*B)/(G*H)]$ is obtained by
quotient forcing over $V[G][H]$. Since $\Hterm*B$ is $\kappa$-c.c., the proof is completed by the fact that any quotient of $\kappa$-c.c.~forcing is
$\kappa$-c.c.\QED

Putting all this together, we now prove Theorem \ref{ForcingNMPccc}.

\Proof Proof of Theorem \ref{ForcingNMPccc}: We assume $V_\delta\elesub V$ and $\kappa$ is a weakly compact cardinal below $\delta$. Let $\cal I$ be
the set of pairs $\<\P*\Qdot,\varphi(\xdot)>$ in $V_\delta$ such that $\P\in V_\kappa$, $\xdot$ is a $\P$-name for an element of $H(\kappa)$ and
$\Qdot\in V_\delta^\P$ is further $\kappa$-c.c.~forcing such that $\one$ forces via $\P*\Qdot$ over $V_\delta$ that $\varphi(\xdot)$ is true in all
$\kappa$-c.c.~forcing extensions of $V_\delta^{\P*\Qdot}$. In this case, let $\Q_{\<\P*\Qdot,\varphi(\xdot)>}$ be the term-forcing for $\Qdot$ over
$\P$, and let $\Q_\infty=\Pi_{\cal I}\Q_{\<\P*\Qdot,\varphi(\xdot)>}$ be the finite support product of these posets. By Lemma
\ref{TermForcingKappaCC}, each factor in this poset is $\kappa$-c.c., and so by Lemma \ref{ProductIsKappaCC} the product $\Q_\infty$ is also
$\kappa$-c.c. Suppose that $G_\infty\of\Q_\infty$ is $V$-generic and consider $V[G_\infty]$. Suppose that $\P$ is some forcing in $V[G_\infty]$
adding an object $\xdot$ in $H(\kappa)$, and that $\varphi(\xdot)$ is forceably necessary for $\kappa$-c.c.~forcing over $V[G_\infty]^\P$. That is,
there is some further $\kappa$-c.c.~forcing $\Qdot$ such that $\varphi(\xdot)$ holds in all $\kappa$-c.c.~extensions of $V[G_\infty]^{\P*\Qdot}$.
Since $\Q_\infty*\P$ adds $\xdot$, there is by Lemma \ref{SmallSuborder} a complete suborder $\P_0\of\Q_\infty*\P$ of size less than $\kappa$ adding
$\xdot$, so let us assume that $\xdot$ is a $\P_0$-name. Since $V[G_\infty]^{\P*\Qdot}$ is a $\kappa$-c.c.~extension of $V^{\P_0}$, it follows that
$\varphi(\xdot)$ is $\kappa$-c.c.~forceably necessary over $V^{\P_0}$, and therefore, by the elementarity $V_\delta^{\P_0}\elesub V^{\P_0}$, it is
$\kappa$-c.c.~forceably necessary over $V_\delta^{\P_0}$. So we may assume that $\Qdot$ was chosen from $V_\delta^{\P_0}$, and that it is forced by
$\P_0*\Qdot$ that $\varphi(\xdot)$ holds in all $\kappa$-c.c.~extensions of $V_\delta^{\P_0*\Qdot}$. Thus, $\<\P_0*\Qdot,\varphi(\xdot)>\in\cal I$,
and $V[G_\infty]$ has a $V$-generic filter $\Hterm$ for the term forcing $\Q_{\<\P_0*\Qdot,\varphi(\xdot)>}$. Since $G_0\in V[G_\infty]$, it follows
by Lemma \ref{TermForcing} that there is a $V[G_0]$-generic filter $H\of\Q=\Qdot_{G_0}$ in $V[G_\infty]$. By the choice of $\Qdot$, we know that
$\varphi(x)$ holds in all $\kappa$-c.c.~extensions of $V_\delta[G_0][H]$, and hence by elementarity in all $\kappa$-c.c.~extensions of $V[G_0][H]$.
Since $V[G_\infty]$ is by Lemma \ref{TermForcingKappaCC} a $\kappa$-c.c.~extension of $V[G_0][H]$, we conclude that $\varphi(x)$ holds there, as
desired. We have established that $V[G_\infty]$ satisfies $\necessary\MP_{\hbox{\!\scriptsize$\kappa$-c.c.}}(H(\kappa))$. It follows that $\kappa$
has become $\omega_1$ in $V[G_\infty]$, because for any $\alpha<\kappa$ the assertion that $\alpha$ is countable is $\kappa$-c.c.~forceably
necessary, and hence true in $V[G_\infty]$. So $\kappa$-c.c.~has become the same as c.c.c., and we have
$V[G_\infty]\satisfies\necessary\MPccc(\R)$.\QED

One may not omit the hypothesis of $V_\delta\elesub V$ in Theorem \ref{ForcingNMPccc}, because if there is a model of $\ZFC$ with a weakly compact
cardinal, then there is such a model having no forcing extension that is a model of $\MPccc$. To see this, following \cite[Corollary
32]{Hamkins2003:MaximalityPrinciple}, suppose $M$ is any model of \ZFC\ with a weakly compact cardinal. Let $N$ be the union of all $L_\theta$ in
this model, where $\theta$ is definable in $L^M$ without parameters. An easy Tarski-Vaught argument shows that $N\elesub L^M$, and so
$N\satisfies\ZFC+V=L$ and the definable ordinals of $N$ are unbounded in $N$. Note that the least weakly compact cardinal of $L^M$ must be in $N$ and
weakly compact there. If $\theta$ is defined by $\varphi(x)$ in $N$, then $\varphi(x)^L$ defines $\theta$ in any forcing extension of $N$.
Consequently, if such an extension $N[G]$ satisfies $\MPccc$, then it would have to satisfy $2^\omega>\theta$, since this is expressible without
parameters and is c.c.c.~forceably necessary. Since the definable ordinals $\theta$ are unbounded in the ordinals of $N$, the value of $2^\omega$ in
$N[G]$ would have to be larger than every ordinal, a contradiction.

This completes the converse direction of the Main Theorem \ref{Main Theorem}. Let us turn now to the forward implication by showing that if the
Necessary Maximality Principle $\necessary\MPccc(\R)$ is consistent with \ZFC, then so is the existence of a weakly compact cardinal. This argument
proceeds to a great measure merely by placing the known results of Theorems \ref{HarringtonShelah} and \ref{NMPcccOmega1Inacc} adjacent to one
another and observing the result.

\Theorem.({Harrington-Shelah \cite[Theorem C.i]{HarringtonShelah1985:SomeExactEquiconsistencyResults}}) If \MA\ holds and $\omega_1$ is inaccessible
to reals, then it is weakly compact in $L$.\label{HarringtonShelah}

\Theorem.({Leibman \cite{Leibman2004:Dissertation}}) If $\necessary\MPccc(\R)$ holds, then $\omega_1$ is inaccessible to
reals.\label{NMPcccOmega1Inacc}

The former result has been widely discussed elsewhere (for example, see \cite{Schindler2000:ForcingAxiomsProjectiveSets}). Leibman's proof of Theorem
\ref{NMPcccOmega1Inacc} proceeds roughly as follows: If $c:\omega\to\omega$ is a Cohen real added over $V$, then by Todor\v cevi\'c's proof of
Shelah's theorem (see \cite[Theorem 28.12]{Jech:SetTheory3rdEdition}) there is in $V[c]$ a Souslin tree $T(c)$ constructed by composing $c$ with each
element of an almost coherent family of injective functions $e_\alpha:\alpha\to\omega$.  If $\omega_1=\omega_1^{L[z]}$ for some real $z$ in $V$, then
there will be such an almost coherent family of injective function in $L[z]$. By using the $L[z]$-least such family, the tree $T(c)$ is seen to be
definable in $V[c]$ from the parameters $z$ and $c$, and furthermore, this definition is absolute to any c.c.c.~extension. The assertion that this
tree has an $\omega_1$ branch, therefore, which uses only the parameters $z$ and $c$, is c.c.c.~forceably necessary (since one can force with the
Souslin tree), but not true in $V[c]$ (since it is a Souslin tree there), contradicting $\necessary\MPccc(\R)$.

Theorems \ref{HarringtonShelah} and \ref{NMPcccOmega1Inacc} now combine to establish the result we need:

\Corollary. If $\necessary\MPccc(\R)$ holds, then $\omega_1$ is weakly compact in $L$.\label{NMPcccOmega1wc}

\Proof: If $\necessary\MPccc(\R)$ holds, then of course it holds in every c.c.c.~extension. Consequently, by Theorem \ref{NMPcccOmega1Inacc}, if
$\necessary\MPccc(\R)$ holds, then $\omega_1$ is inaccessible to reals in every c.c.c.~extension. Since there is such an extension satisfying \MA, it
follows by applying Theorem \ref{HarringtonShelah} in that extension that $\omega_1$ is weakly compact in $L$, as desired.\QED

Let us give a second proof, along a different route. The first step is an easy induction on formulas:

\Theorem. If $\necessary\MPccc(\R)$ holds, then projective assertions are absolute by c.c.c.~forcing.

\Proof: The method of \cite[Theorem 19]{Hamkins2003:MaximalityPrinciple} works generally. First notice that if $\necessary\MPccc(\R)$ holds, then it
holds in all c.c.c.~extensions. Suppose inductively that a projective assertion $\varphi(x,y)$ is absolute from any c.c.c.~extension to any further
c.c.c.~extension. Boolean combinations are easily preserved as well, so it suffices to consider the existential case. If $\exists x\,\varphi(x,a)$ is
true in $V^\P$, then by substituting the witness into place, this is preserved by induction to any further extension $V^{\P*\Qdot}$. Conversely, if
$V^{\P*\Qdot}\satisfies\exists x\,\varphi(x,a)$, then the existence of a witness is forceably necessary over $V^\P$, and hence true in $V^\P$, as
desired.\QED

Since c.c.c.~projective absoluteness is known to be equiconsistent with the existence of a weakly compact cardinal, Corollary \ref{NMPcccOmega1Inacc}
now follows. \cite{BagariaFriedman2001:GenericAbsoluteness} shows already that c.c.c.~$\Sigma^1_4$ absoluteness is equiconsistent with a weakly
compact cardinal.

Let us close the paper with an application of our method to the case of the original Maximality Principle \MP, without restricting to c.c.c.~forcing.
By \cite[Theorem 12]{Hamkins2003:MaximalityPrinciple} we know that it is consistent that $\MP(\R)$ is indestructible by the forcing to add Cohen
reals, but in fact $\MP(\R)$ is always indestructible.

\Theorem. The boldface Maximality Principle $\MP(\R)$, if true, is indestructible by the forcing to add any number of Cohen
reals.\label{MPRCohenIndestructible}

\Proof: Suppose first that $\MP(\R)$ holds in $V$ and $V[c]$ is a generic extension obtained by adding a Cohen real $c$. To show that $V[c]$ models
$\MP(\R)$, suppose that $x$ is a real in $V[c]$ and $\varphi(x)$ is forceably necessary in $V[c]$. Thus, there is a poset $\Q$ in $V[c]$ such that if
$G\of\Q$ is $V[c]$-generic, then $\varphi(x)$ holds in $V[c][G]$ and all extensions. Let $\dot x$ be a name for $x$ and $\Qdot$ a name for $\Q$, such
that there is a condition $p_0$ in $c$ forcing that $\Qdot$ makes $\varphi(\dot x)$ necessary. Let $\Qterm$ be the term forcing poset of $\Qdot$ over
the Cohen real forcing $\C$, and suppose $\Gterm\of\Qterm$ is $V[c]$-generic. The model $V[c][\Gterm]=V[\Gterm][c]$ has a $V$-generic $c\of\C$ and
$\Gterm\of\Qterm$, so by the fundamental property of term forcing Lemma \ref{TermForcing}, there is a $V[c]$-generic filter $G\of\Q$ in
$V[c][\Gterm]$ and we may view the extension as $V[\Gterm][c]=V[c][G][\Gterm/G]$, that is, as an extension of $V[c][G]$. Since $\varphi(x)$ was made
necessary in $V[c][G]$, it follows that $V[\Gterm][c]\satisfies\varphi(x)$. In fact we know that $p_0\forces_\C\varphi(\dot x)$ in $V[\Gterm]$.
Furthermore, if $V[\Gterm][H]$ is some other forcing extension of $V[\Gterm]$, and $c$ is any $V[\Gterm][H]$ generic below $p_0$, then we may
rearrange the resulting extension as $V[\Gterm][H][c]=V[\Gterm][c][H]=V[c][G][\Gterm/G][H]$, which still satisfies $\varphi(x)$ since it is an
extension of $V[c][G]$. Therefore, $p_0\forces_\C\varphi(\dot x)$ in $V[\Gterm][H]$ as well. In short, we have proved that the assertion
``$p_0\forces_\C\varphi(\dot x)$'' is necessary in $V[\Gterm]$, and therefore forceably necessary in $V$. Since the parameters $\C$, $p_0$ and $\dot
x$ in this assertion are all hereditarily countable in $V$, we conclude by $\MP(\R)$ in $V$ that the assertion must be true in $V$. Thus, since $p_0$
is in $c$, we conclude that $V[c]\satisfies\varphi(x)$, without any need for the term forcing, as desired.

Now we prove the full result of the theorem. Suppose that $G\of\add(\omega,\kappa)$ is $V$-generic for the forcing to add $\kappa$ many Cohen reals.
By the countable chain condition, every real $x$ of $V[G]$ is in $V[G\intersect a]$ for some countable set $a\in V$. Since $G\intersect a$ is
isomorphic to adding a single Cohen real, we know by the previous paragraph that $\MP(\R^{V[G\intersect a]})$ is true in $V[G\intersect a]$. In
particular, any instance of the \MP\ scheme involving the parameter $x$ holds in $V[G\intersect a]$. Since $V[G]$ is a forcing extension of
$V[G\intersect a]$, it follows that if $\varphi(x)$ is forceably necessary over $V[G]$, then it is forceably necessary over $V[G\intersect a]$, and
hence necessary in $V[G\intersect a]$, and hence still necessary in $V[G]$. So we have established every instance of the $\MP(\R)$ scheme in $V[G]$,
as desired.\QED

The theorem applies more generally, of course, to any forcing extension all of whose reals are captured by Cohen forcing. We note the contrast of
Theorem \ref{MPRCohenIndestructible} with the consequence of Leibman's proof of Observation \ref{GeorgeObservation}, that $\MPccc(H(2^\omega))$ is
always destroyed by the forcing to add a Cohen real. It follows that if $\necessary\MPccc(\R)$ is consistent, it is consistent with the failure of
$\MPccc(H(2^\omega))$.

\bibliographystyle{alpha}
\bibliography{MathBiblio,HamkinsBiblio}

\end{document}